\documentclass[12pt,reqno]{amsart}
\usepackage{latexsym,amsmath,amsfonts,amssymb,amsthm,color}
\theoremstyle{plain}
\newtheorem{Thm}{Theorem}
\newtheorem{Lem}{Lemma}
\newtheorem{Cor}{Corollary}

\theoremstyle{definition}

\theoremstyle{remark}

\def\N{\mathbb N}

\def\C{\mathbb C}

\def\1{{\bf 1}}

\def\M#1#2#3#4{\left[\begin{matrix} #1&#2\\#3&#4\end{matrix}\right]_n}

\begin{document}
\hbox{Final version of the paper (GCOM-D-09-00021)
accepted by Graphs and Combin.}
\medskip
\title{Symmetric identities for Euler polynomials}
\author{Yong Zhang, Zhi-Wei Sun and Hao Pan}
\address{Department of Mathematics, Nanjing University, Nanjing 210093,
People's Republic of China
\newline \indent
 {\it E-mail addresses}: {\tt
yongzhang1982@163.com, zwsun@nju.edu.cn, haopan79@yahoo.com.cn}}

\thanks {The second author is responsible for communications, and supported by
the National Natural Science Foundation (grant 10871087)
and the Overseas Cooperation Fund (grant 10928101) of China.}
 \maketitle

\vskip 10pt

\begin{abstract} In this paper we establish two symmetric identities on sums of products of Euler
polynomials.
\end{abstract}

\section{Introduction}
\setcounter{Lem}{0}\setcounter{Thm}{0}\setcounter{Cor}{0}
\setcounter{equation}{0}

 The Bernoulli numbers $B_0,B_1,B_2,\ldots$ are rational numbers given by
 $$B_0=1,\ \mbox{and}\ \sum_{k=0}^n\binom{n+1}kB_k=0\ \ \mbox{for}\ n=1,2,3,\ldots.$$
The Euler numbers $E_0,E_1,E_2,\ldots$ are integers determined by
$$E_0=1, \ \mbox{and}\ \sum^n_{\substack {k=0\\2\mid n-k}}\binom{n}kE_k=0\ \ \mbox{for}\ n=1,2,3,\ldots.$$

Let $\N=\{0,1,2,\ldots\}$. The Bernoulli polynomials $B_n(x)\
(n\in\N)$ and the Euler polynomials $E_n(x)\ (n\in\N)$ are defined
by
$$B_n(x)=\sum_{k=0}^n\binom{n}{k}B_kx^{n-k} \ \mbox{and}\
E_n(x)=\sum_{k=0}^n\binom{n}{k}\frac{E_k}{2^k}\left(x-\frac{1}{2}\right)^{n-k}.$$
It is well known that
$$\Delta(B_n(x))=nx^{n-1}\ \mbox{and}\ \Delta^*(E_n(x))=2x^n$$
for all $n\in\N$, where we set
$$\Delta(P(x))=P(x+1)-P(x)\ \mbox{and}\ \Delta^*(P(x))=P(x+1)+P(x)$$
for any polynomial $P(x)$. Bernoulli and Euler numbers and
polynomials play important roles in many fields including number
theory and combinatorics.

 In 2006 Z. W. Sun and H. Pan \cite{SunPan06} established the following theorem which unifies
 many  curious identities
 concerning Bernoulli and Euler numbers and polynomials.

\begin{Thm} [Sun and Pan, 2006] Let $n$ be a positive integer and
let $x+y+z=1$.

{\rm (i)} If $r,s,t$ are complex numbers with $r+s+t=n$, then we
have the symmetric relation
$$r\M stxy+s\M tryz+t\M rszx=0$$
where
$$\M stxy:=\sum_{k=0}^n(-1)^k\binom sk\binom t{n-k}B_{n-k}(x)B_k(y).$$

{\rm (ii)} If $r+s+t=n-1$, then
\begin{eqnarray*}
\label{ab}
\lefteqn{\frac{r}{2}\sum_{l=0}^{n-1}(-1)^l\binom{s}{l}\binom{t}{n-1-l}E_l(y)E_{n-1-l}(x)}\notag\\
&=&\sum_{k=0}^n(-1)^k\binom{r}{k}\binom{s}{n-k}B_k(x)E_{n-k}(z)\notag\\
&&-(-1)^n\sum_{k=0}^n(-1)^k\binom{r}{k}\binom{t}{n-k}B_k(y)E_{n-k}(z).
\end{eqnarray*}
\end{Thm}

Recently, by a sophisticated application of the generating function
method, A. M. Fu, H. Pan and F. Zhang \cite{FPZ} extended Theorem
1.1(i) of Sun and Pan to an identity on sums of products of $m\ge2$
Bernoulli polynomials.

In this paper we obtain a general identity only involving Euler
polynomials and also give an extension of Theorem 1.1(ii) which
involves both Bernoulli and Euler polynomials.

\begin{Thm} \label{t1} Let $m$ and $n$ be positive
integers, and let $r_0,r_1,\ldots,r_m$ be complex numbers with
$r_0+r_1+\cdots+r_m=n-1$.

{\rm(i)} If $m$ is odd, then we have the symmetric relation
\begin{align}
\label{euler}
&\sum_{\substack{k_1,\cdots,k_m\geq 0\\
k_1+\cdots+k_m=n}}\prod_{j=1}^m\binom{r_j}{k_j}E_{k_j}(x_j)\notag\\
=&-\sum_{i=1}^m(-1)^i\sum_{\substack{k_1,\ldots,k_m\geq 0\\
k_1+\cdots+k_m=n}}\binom{r_0}{k_i}E_{k_i}(1-x_i)\prod_{\substack{1\leq
j\leq m\\j\not=i}}\binom{r_j}{k_j}E_{k_j}(x_j-x_i+\1_{j>i}),
\end{align}
where $\1_{j>i}$ takes $1$ or $0$ according as $j>i$ or not.

{\rm(ii)} If $m$ is even, then
\begin{align}
\label{bernoulli}
&\frac{r_0}{2}\sum_{\substack{k_1,\cdots,k_m\geq 0\\
k_1+\cdots+k_m=n-1}}\prod_{j=1}^m\binom{r_j}{k_j}E_{k_j}(x_j)\notag\\
=&\sum_{i=1}^m(-1)^i\sum_{\substack{k_1,\ldots,k_m\geq 0\\
k_1+\cdots+k_m=n}}\binom{r_0}{k_i}B_{k_i}(1-x_i)\prod_{\substack{1\leq
j\leq m\\j\not=i}}\binom{r_j}{k_j}E_{k_j}(x_j-x_i+\1_{j>i}).
\end{align}
\end{Thm}

\medskip
{\it Remark}\ 1.1. If $r+s+t=n-1$, then  (\ref{bernoulli}) in the
case $m=2$ gives
\begin{align*}
&\frac{r}{2}\sum_{k=0}^{n-1}\binom{s}{k}E_k(1-y)\binom{t}{n-1-k}E_{n-1-k}(x)\\
=&-\sum_{k=0}^n\binom{r}{k}B_k(1-(1-y))\binom{t}{n-k}E_{n-k}(x-(1-y)+1)\\
&+\sum_{k=0}^n\binom{r}{k}B_k(1-x)\binom{s}{n-k}E_{n-k}((1-y)-x)\\
=&-(-1)^n\sum_{k=0}^n(-1)^k\binom{t}{n-k}E_{n-k}(1-x-y)\binom{r}{k}B_{k}(y)\\
&+\sum_{k=0}^n(-1)^k\binom{r}{k}B_k(x)\binom{s}{n-k}E_{n-k}(1-x-y),
\end{align*}
which is equivalent to the identity of Sun and Pan in Theorem
1.1(ii) since $E_k(1-x)=(-1)^kE_k(x)$.
\medskip

 Our proof of Theorem \ref{t1} given in the next section involves the difference operator
 $\Delta$ and its companion operator $\Delta^*$. We can also show
 Theorem \ref{t1} via the generating function approach.

Let $k$ be any nonnegative integer. It is well known that  $B_k=0$ if $k$ is odd and greater than one.
By \cite[pp. 804-808]{AS},
$$B_k\left(\frac12\right)=(2^{1-k}-1)B_k\ \ \mbox{and}\ \ E_k(x)
=\frac2{k+1}\left(B_{k+1}(x)-2^{k+1}B_{k+1}\left(\frac x2\right)\right).$$
 Thus
$$(-1)^kE_k(1)=E_k(0)=2(1-2^{k+1})\frac{B_{k+1}}{k+1}.$$
In view of these, Theorem \ref{t1} in the case $x_1=\cdots=x_m=1/2$ yields the following consequence
involving Euler numbers and Bernoulli numbers.

\begin{Cor} \label{c1} Let $m$ and $n$ be positive
integers, and let $r_0,r_1,\ldots,r_m$ be complex numbers with
$r_0+r_1+\cdots+r_m=n-1$.

{\rm(i)} If $m$ is odd, then
\begin{align}&(-1)^n\sum_{\substack{k_1,\cdots,k_m\geq 0\\
k_1+\cdots+k_m=n}}\prod_{j=1}^m\binom{r_j}{k_j}E_{k_j}
\notag\\=&\sum_{i=1}^m(-1)^i\sum_{\substack{k_1,\ldots,k_m\geq 0\\
k_1+\cdots+k_m=n}}(-1)^{|\{i<j\le m:\ k_j>0\}|}\binom{r_0}{k_i}E_{k_i}\prod_{\substack{1\leq j\leq
m\\j\not=i}}\binom{r_j}{k_j}\tilde B_{k_j+1},
\end{align}
where $\tilde B_k=2^k(2^k-1)B_k/k$ for $k=1,2,3,\ldots$.

{\rm(ii)} If $m$ is even, then
\begin{align}
&(-1)^nr_0\sum_{\substack{k_1,\cdots,k_m\geq 0\\
k_1+\cdots+k_m=n-1}}\prod_{j=1}^m\binom{r_j}{k_j}E_{k_j}
\notag\\=&\sum_{i=1}^m(-1)^i\sum_{\substack{k_1,\ldots,k_m\geq 0\\
k_1+\cdots+k_m=n}}(-1)^{|\{i<j\le m:\ k_j>0\}|}\binom{r_0}{k_i}(2^{k_i}-2)B_{k_i}
\prod_{\substack{1\leq j\leq m\\j\not=i}}\binom{r_j}{k_j}\tilde B_{k_j+1}.
\end{align}
\end{Cor}

\section{Proof of Theorem \ref{t1}}
\setcounter{Lem}{0}\setcounter{Thm}{0}\setcounter{Cor}{0}
\setcounter{equation}{0}

As usual we let $\C$ denote the field of complex numbers.
By \cite[Lemma 3.1]{sunpan05}, for $P(x),Q(x)\in\C[x]$, we have
$P(x)=Q(x)$ if $\Delta^*(P(x))=\Delta^*(Q(x))$. This property will
play a central role in our proof of Theorem \ref{t1}.

\begin{Lem}
\label{bL1}  Let $P_{1}(x),\cdots,P_{m}(x)\in\C[x]$.
 Then
\begin{align*}&P_1(x)\sum_{1<i\le
m}(-1)^i\Delta^*(P_i(x))\prod_{\substack{1<j\leq
m\\j\not=i}}P_{j}(x+\1_{j<i})
\\=&\begin{cases}\Delta^*(P_1(x)\cdots
P_m(x))-\Delta^*(P_1(x))P_2(x+1)\cdots P_m(x+1)&\mbox{if}\ 2\nmid
m,\\\Delta^*(P_1(x)\cdots P_m(x))-\Delta(P_1(x))P_2(x+1)\cdots
P_m(x+1)&\mbox{if}\ 2\mid m.\end{cases}
\end{align*}
\end{Lem}

\begin{proof} Observe that
$$\aligned&\sum_{1<i\le m}(-1)^i\Delta^*(P_i(x))
\prod_{\substack{1<j\leq m\\j\not=i}}P_{j}(x+\1_{j<i})
\\=&\sum_{1<i\le m}\bigg((-1)^i\prod_{1<j\leq m}P_{j}(x+\1_{j<i})
-(-1)^{i+1}\prod_{1<j\leq m}P_{j}(x+\1_{j<i+1})\bigg)
\\=&(-1)^2\prod_{1<j\le m}P_j(x)-(-1)^{m+1}\prod_{1<j\le m}P_j(x+1).
\endaligned$$
Therefore
$$\aligned&P_1(x)\sum_{1<i\le m}(-1)^i\Delta^*(P_i(x))
\prod_{\substack{1<j\leq m\\j\not=i}}P_{j}(x+\1_{j<i})
\\=&P_1(x)\cdots P_m(x)+(-1)^mP_1(x)\prod_{1<j\le m}P_j(x)
\\=&\Delta^*(P_1(x)\cdots P_m(x))-(P_1(x+1)+(-1)^{m-1}P_1(x))\prod_{1<j\le
m}P_j(x).
\endaligned$$
This proves the desired identity.
\end{proof}

\begin{Lem}
\label{cL1}
Let $a_0,\bar a_0,a_1,\bar a_1\ldots,a_n,\bar a_n$ be complex numbers, and set
$$A_k(t)=\sum_{l=0}^k\binom{k}{l}(-1)^la_lt^{k-l}\ \ \mbox{and}
\ \ \bar A_k(t)=\sum_{l=0}^k\binom{k}{l}(-1)^l\bar a_lt^{k-l}$$
for $k=0,\ldots,n$.
Let $r_0+r_1+\cdots+r_m=n-1$. Then
\begin{align}
\label{cczhang}
&\sum_{\substack{k_1,\ldots,k_m\geq0\\k_1+\cdots+k_m=n}}\binom{r_0}{k_1}(-x_1)^{k_1}
\prod_{j=2}^m\binom{r_j}{k_j}A_{k_j}(x_j-x_1)\notag\\
=&\sum_{\substack{k_1,\ldots,k_m\geq0\\k_1+\cdots+k_m=n}}\binom{r_1}{k_1}{x_1}^{k_1}
\prod_{j=2}^m\binom{r_j}{k_j}A_{k_j}(x_j).
\end{align}
Also, for any $i=2,\ldots,m$ we have
\begin{align}
\label{ccczhang}
&\sum_{\substack{k_1,\ldots,k_m\geq0\\k_1+\cdots+k_m=n}}\binom{r_0}{k_1}
A_{k_1}(-x_1)\binom{r_i}{k_i}(x_i-x_1)^{k_i}\prod_{\substack{2\leq
j\leq m\\ j\not=i}}\binom{r_j}{k_j}\bar A_{k_j}(x_j-x_1)\notag\\
=&\sum_{\substack{k_1,\ldots,k_m\geq0\\k_1+\cdots+k_m=n}}\binom{r_1}{k_1}(x_1-x_i)^{k_1}\binom{r_0}{k_i}
A_{k_i}(-x_i)\prod_{\substack{2\leq j\leq m\\
j\not=i}}\binom{r_j}{k_j}\bar A_{k_j}(x_j-x_i).
\end{align}
\end{Lem}

\begin{proof}
By Remark 1.1 of Sun \cite{sun},
$$A_k(x+y)=\sum_{l=0}^k\binom{k}{l}x^{k-l}A_l(y)
\ \ \mbox{and}\ \ \bar A_k(x+y)=\sum_{l=0}^k\binom{k}{l}x^{k-l}\bar A_l(y)$$
for every $k=0,\ldots,n$.
Observe that
\begin{align*}
&\sum_{\substack{k_1,\ldots,k_m\geq0\\k_1+\cdots+k_m=n}}\binom{r_0}{k_1}(-x_1)^{k_1}
\prod_{j=2}^m\binom{r_j}{k_j}A_{k_j}(x_j-x_1)\notag\\
=&\sum_{\substack{k_1,\ldots,k_m\geq0\\k_1+\cdots+k_m=n}}\binom{r_0}{k_1}(-x_1)^{k_1}
\prod_{j=2}^m\binom{r_j}{k_j}\sum_{l_j=0}^{k_j}\binom{k_j}{l_j}(-x_1)^{k_j-l_j}A_{l_j}(x_j)\notag\\
=&\sum_{\substack{l_1,\ldots,l_m\geq0\\l_1+\cdots+l_m=n}}(-x_1)^{l_1}\prod_{j=2}^m
\binom{r_j}{l_j}A_{l_j}(x_j) \sum_{\substack{k_1\geq0,\ k_j\geq l_j\
(1<j\leq m)\\k_1+\cdots+k_m=n}}\binom{r_0}{k_1}
\prod_{j=2}^m\binom{r_j-l_j}{k_j-l_j}.\notag
\end{align*}
Given $l_1,\ldots,l_m\in\N$ with $l_1+\cdots+l_m=n$,
by the Chu-Vandermonde convolution identity (cf. \cite[(5.22)]{GKP}), we have
\begin{align*}
&\sum_{\substack{k_1\geq0,\,k_j\geq l_j\ (1<j\leq
m)\\k_1+\cdots+k_m=n}}\binom{r_0}{k_1}\prod_{j=2}^m\binom{r_j-l_j}{k_j-l_j}\\
=&\binom{r_0+(r_2-l_2)+\cdots+(r_m-l_m)}{n-l_2-\cdots-l_m}
=\binom{l_1-1-r_1}{l_1}=(-1)^{l_1}\binom{r_1}{l_1}.
\end{align*}
So (\ref{cczhang}) follows.

(\ref{ccczhang}) can be proved similarly. Let $\Sigma$ denote the left-hand side of (\ref{ccczhang}). Then
\begin{align*}
\Sigma
=&\sum_{\substack{k_1,\ldots,k_m\geq0\\k_1+\cdots+k_m=n}}
\binom{r_0}{k_1}\sum_{l_i=0}^{k_1}\binom{k_1}{l_i}(x_i-x_1)^{k_1-l_i}
A_{l_i}(-x_i)\binom{r_i}{k_i}(x_i-x_1)^{k_i}\notag\\
&\times\prod_{\substack{1<j\leq m\\j\not=i}}\binom{r_j}{k_j}\sum_{l_j=0}^{k_j}\binom{k_j}{l_j}(x_i-x_1)^{k_j-l_j}
\bar A_{l_j}(x_j-x_i)\notag\\
=&\sum_{\substack{l_1,\ldots,l_m\geq0\\l_1+\cdots+l_m=n}}(x_i-x_1)^{l_1}\binom{r_0}{l_i}A_{l_i}(-x_i)
\prod_{\substack{1<j \leq m\\j\not=i}}\binom{r_j}{l_j}\bar A_{l_j}(x_j-x_i)\notag\\
&\times\sum_{\substack{k_j\geq l_j\,(1\leq j \leq m\ \&\
j\not=i)\\k_i\geq0,\ k_1+\cdots+k_m=n}}\binom{r_0-l_i}{k_1-l_i}
\binom{r_i}{k_i}\prod_{\substack{1<j\leq m\\j\not=i}}\binom{r_j-l_j}{k_j-l_j}\notag\\
=&\sum_{\substack{l_1,\ldots,l_m\geq0\\l_1+\cdots+l_m=n}}(x_i-x_1)^{l_1}\binom{r_0}{l_i}A_{l_i}(-x_i)
\prod_{\substack{1<j \leq
m\\j\not=i}}\binom{r_j}{l_j}\bar A_{l_j}(x_j-x_i)\times
(-1)^{l_1}\binom{r_1}{l_1}.
\end{align*}
This concludes the proof.
\end{proof}

\medskip
\noindent{\it Remark}\ 2.1. If we set $a_l=(-1)^lB_l$ and $\bar a_l=(-1)^lE_l(0)$ for $l=0,\ldots,n$
in Lemma 2.2, then $A_k(t)=B_k(t)$ and $\bar A_k(t)=E_k(t)$ for any $k=0,\ldots,n$.
\medskip

\begin{proof}[Proof of Theorem \ref{t1}]
We fix $x_{2},\ldots,x_{m}$.

{\rm(i)} Suppose that $m$ is odd. Set
$$P(x_{1})=\sum_{\substack{k_1,\ldots,k_m\geqq0\\
k_1+\cdots+k_m=n}}\binom{r_0}{k_1}E_{k_1}(1-x_1)\prod_{j=2}^m\binom{r_j}{k_j}E_{k_j}(x_j-x_1+1).$$
Applying Lemma \ref{bL1}, we get
\begin{align*}
&\Delta^*(P(x_1))\notag\\=&\sum_{i=2}^m(-1)^i\sum_{\substack{k_1,\ldots,k_m\geqq0
\\k_1+\cdots+k_m=n}}\binom{r_0}{k_1}E_{k_1}(1-x_1)\binom{r_i}{k_i}2(x_i-x_1)^{k_i}\prod_{\substack{2\leq
j\leq m\\j\not=i}}\binom{r_j}{k_j}E_{k_j}(x_j-x_1+\1_{j>i})\notag\\
&+\sum_{\substack{k_1,\ldots,k_m\geqq0\\k_1+\cdots+k_m=n}}\binom{r_0}{k_1}2(-x_1)^{k_1}\prod_{j=2}^m\binom{r_j}{k_j}
E_{k_j}(x_j-x_1).
\end{align*}
With the help of Lemma \ref{cL1}, we have
\begin{align*}
&\Delta^*(P(x_1))\notag\\
=&2\sum_{i=2}^m(-1)^i\sum_{\substack{k_1,\ldots,k_m\geqq0\\k_1+\cdots+k_m=n}}
\binom{r_1}{k_1}(x_1-x_i)^{k_1}\binom{r_0}{k_i}E_{k_i}(1-x_i)
\prod_{\substack{2\leq j\leq
m\\j\not=i}}\binom{r_j}{k_j}E_{k_j}(x_j-x_i+\1_{j>i})\notag\\
&+2\sum_{\substack{k_1,\ldots,k_m\geqq0\\k_1+\cdots+k_m=n}}\binom{r_1}{k_1}x_1^{k_1}\prod_{j=2}^m
\binom{r_j}{k_j}E_{k_j}(x_j).
\end{align*}
It follows that $\Delta^*(P(x_1))=\Delta^*(Q(x_1))$, where
\begin{align*}
Q(x_1)=&\sum_{1<i\le m}(-1)^i\sum_{\substack{k_1,\ldots,k_m\geq 0\\
k_1+\cdots+k_m=n}}\binom{r_0}{k_i}E_{k_i}(1-x_i)\prod_{\substack{1\leq
j\leq m\\j\not=i}}\binom{r_j}{k_j}E_{k_j}(x_j-x_i+\1_{j>i})\notag\\
&+\sum_{\substack{k_1,\cdots,k_m\geq 0\\
k_1+\cdots+k_m=n}}\prod_{j=1}^m\binom{r_j}{k_j}E_{k_j}(x_j).
\end{align*}
Therefore $P(x_1)=Q(x_1)$ by \cite[Lemma 3.1]{sunpan05}. This proves
(\ref{euler}).
\medskip

{\rm(ii)} Now assume that $m$ is even. Define
$$P(x_{1})=\sum_{\substack{k_1,\ldots,k_m\geqq0\\
k_1+\cdots+k_m=n}}\binom{r_0}{k_1}B_{k_1}(1-x_1)\prod_{j=2}^m\binom{r_j}{k_j}E_{k_j}(x_j-x_1+1).
$$
For $k_1=0,1,2,\ldots$, clearly
$$\binom{r_0}{k_1}(B_{k_1}(1-(x_1+1))-B_{k_1}(1-x_1))=-\binom{r_0}{k_1}k_1(-x_1)^{k_1-1}
=-r_0\binom{r_0-1}{k_1-1}(-x_1)^{k_1-1}.$$
(As usual $\binom x{-1}$ is regarded as $0$.)
Thus, by Lemma \ref{bL1} we have
\begin{align*}
&\Delta^*(P(x_1))\notag\\=&2\sum_{i=2}^m(-1)^i\sum_{\substack{k_1,\ldots,k_m\ge0
\\k_1+\cdots+k_m=n}}\binom{r_0}{k_1}B_{k_1}(1-x_1)\binom{r_i}{k_i}(x_i-x_1)^{k_i}\prod_{\substack{1<
j\leq m\\j\not=i}}\binom{r_j}{k_j}E_{k_j}(x_j-x_1+\1_{j>i})\notag\\
&-r_0\sum_{\substack{k_1,\ldots,k_m\ge0\\k_1+\cdots+k_m=n-1}}\binom{r_0-1}{k_1}(-x_1)^{k_1}\prod_{1<j\le
m}\binom{r_j}{k_j} E_{k_j}(x_j-x_1).
\end{align*}
With the help of Lemma \ref{cL1},
\begin{align*}
&\Delta^*(P(x_1))\notag\\
=&2\sum_{i=2}^m(-1)^i\sum_{\substack{k_1,\ldots,k_m\ge0\\k_1+\cdots+k_m=n}}
\binom{r_1}{k_1}(x_1-x_i)^{k_1}\binom{r_0}{k_i}B_{k_i}(1-x_i)
\prod_{\substack{1<j\leq
m\\j\not=i}}\binom{r_j}{k_j}E_{k_j}(x_j-x_i+\1_{j>i})\notag\\
&-r_0\sum_{\substack{k_1,\ldots,k_m\ge0\\k_1+\cdots+k_m=n-1}}\binom{r_1}{k_1}x_1^{k_1}\prod_{1<j\le
m}\binom{r_j}{k_j}E_{k_j}(x_j).
\end{align*}
So we have $\Delta^*(P(x_1))=\Delta^*(Q(x_1))$, where
\begin{align*}
Q(x_1)=&\sum_{i=2}^m(-1)^i\sum_{\substack{k_1,\ldots,k_m\geq 0\\
k_1+\cdots+k_m=n}}\binom{r_0}{k_i}B_{k_i}(1-x_i)\prod_{\substack{1\leq
j\leq m\\j\not=i}}\binom{r_j}{k_j}E_{k_j}(x_j-x_i+\1_{j>i})\notag\\
&-\frac{r_0}{2}\sum_{\substack{k_1,\cdots,k_m\geq 0\\
k_1+\cdots+k_m=n-1}}\prod_{j=1}^m\binom{r_j}{k_j}E_{k_j}(x_j).
\end{align*}
Therefore, $P(x_1)$ coincides with $Q(x_1)$ by \cite[Lemma
3.1]{sunpan05}. So (\ref{bernoulli}) holds. This concludes the proof.
\end{proof}

\medskip
\noindent{\bf Acknowledgment}. We thank the two referees for their helpful comments.

\end{document}